\documentclass{gtmon_a}
\pdfoutput=1

\usepackage[all]{xy}
\usepackage{graphicx}

\proceedingstitle{Proceedings of the School and Conference in Algebraic
Topology (The Vietnam National University, Hanoi, 9-20 August 2004)}
\conferencestart{9 August 2004}
\conferenceend{20 August 2004}
\conferencename{School and Conference in Algebraic Topology}
\conferencelocation{Vietnam National University, Hanoi, Vietnam}

\editor{John Hubbuck}
\givenname{John}
\surname{Hubbuck}

\editor{Nguy\~\ecircumflex{}n H V H\uhorn{}ng}
\givenname{H\uhorn{}ng}
\surname{Nguy\~\ecircumflex{}n}

\editor{Lionel Schwartz}
\givenname{Lionel}
\surname{Schwartz}

\title{Lectures on the stable homotopy of $BG$}

\author{Stewart Priddy}
\givenname{Stewart}
\surname{Priddy}
\address{Department of Mathematics\\
Northwestern University\\\newline
Evanston, IL 60208\\USA}
\email{priddy@math.northwestern.edu}
\urladdr{}

\volumenumber{11}
\issuenumber{}
\publicationyear{2007}
\papernumber{13}
\startpage{289}
\endpage{308}

\doi{}
\MR{}
\Zbl{}

\arxivreference{0710.5942}

\keyword{classifying space}
\keyword{stable splittings}
\keyword{transfer}
\subject{primary}{msc2000}{55P42}
\subject{secondary}{msc2000}{55R35}
\subject{secondary}{msc2000}{20C20}

\received{17 January 2005}
\revised{20 June 2005}
\accepted{20 June 2005}
\published{14 November 2007}
\publishedonline{14 November 2007}
\proposed{}
\seconded{}
\corresponding{}
\editor{}
\version{}

\let\xysavmatrix\xymatrix
\def\xymatrix{\disablesubscriptcorrection\xysavmatrix}
\AtBeginDocument{\let\bar\wbar\let\tilde\wtilde\let\hat\what}
\makeop{tr}
\makeop{Id}
\makeop{End}
\makeop{Cen}
\makeop{St}
\makeop{Reg}

\theoremstyle{plain}
\newtheorem{thm}{Theorem}

\newtheorem{lem}{Lemma}
\newtheorem{cor}{Corollary}
\newtheorem{prop}{Proposition}
\theoremstyle{definition}
\newtheorem*{defn}{Definition}
\newtheorem{exm}{Example}

\makeautorefname{exm}{Example}
\DeclareMathOperator{\Aut}{Aut}
\DeclareMathOperator{\Inj}{Inj}
\DeclareMathOperator{\Hom}{Hom}
\DeclareMathOperator{\Out}{Out}
\DeclareMathOperator{\Rep}{Rep}
\DeclareMathOperator{\colim}{colim}
\DeclareMathOperator{\F}{\mathbb F}

\begin{document}

\begin{asciiabstract}
 This paper is a survey of the stable homotopy theory of BG for
G a finite group. It is based on a series of lectures given at
the Summer School associated with the Topology Conference at the Vietnam National University, Hanoi, August 2004.
\end{asciiabstract}

\begin{abstract}
This paper is a survey of the stable homotopy theory of $BG$ for
$G$ a finite group. It is based on a series of lectures given at
the Summer School associated with the Topology Conference at the Vietnam National University, Hanoi, August 2004.
\end{abstract}

\maketitle

Let $G$ be a finite group. Our goal is to study the stable homotopy of the classifying space $BG$ completed at
some prime $p$. For ease of notation, we shall always assume that any space in question has been $p$--completed.
Our fundamental approach is to decompose the stable type of $BG$ into its various summands. This is useful in
addressing many questions in homotopy theory especially when the summands can be identified with simpler or at
least better known spaces or spectra. It turns out that the summands of $BG$ appear at various levels related to
the subgroup lattice of $G$. Moreover since we are working at a prime, the modular representation theory of
automorphism groups of $p$--subgroups of $G$ plays a key role. These automorphisms arise from the normalizers of
these subgroups exactly as they do in $p$--local group theory. The end result is that a complete stable
decomposition of $BG$ into indecomposable summands can be described (\fullref{t6}) and its stable homotopy type can be
characterized algebraically (\fullref{t7}) in terms of simple modules of automorphism groups.

This paper is a slightly expanded version of lectures given at the International School of the Hanoi Conference on
Algebraic Topology, August 2004. The author is extremely grateful to Mike Hill for taking notes and producing a
TeX document which formed the basis of the present work. Additional comments and references have been added to
make the exposition reasonably self-contained. Many of the results of this paper were obtained jointly with my
longtime collaborator, John Martino.

\section{Preliminaries}
\subsection{What do we mean by ``stable homotopy''}
Given a pointed space $(X,*),$ let $\Sigma X=X\times
I/X\times\{0,1\}\cup\{*\}\times I$ denoted the reduced suspension.
We can represent this pictorially as:
\begin{figure}[ht!]
\centering
\includegraphics[height=1in]{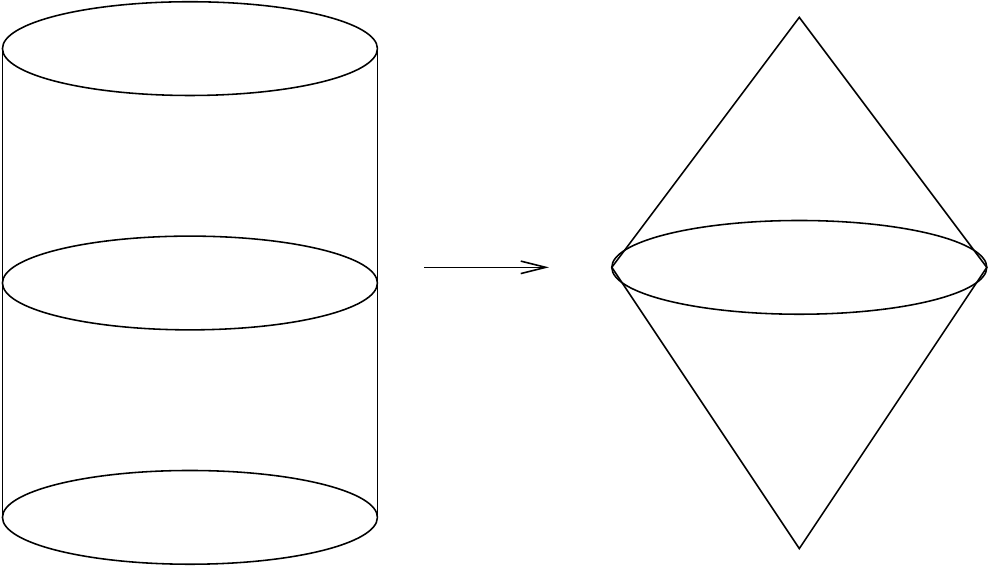}
\end{figure}

One can also quickly check that $\Sigma S^n=S^{n+1}$, where $S^n$ is the $n$--sphere.

Now, $\tilde{H}_*(\Sigma X)=\tilde{H}_{*-1}(X)$, and for any space $Y$, $[\Sigma X,Y]$ is a group. The composition
is as defined for homotopy groups: if we have two homotopy classes $[f],[g]$, then we let $[f]\cdot[g]$ denote the
composite $[\nabla\circ (f,g)\circ\pi']$, where $\pi'\co \Sigma X\to \Sigma X\vee\Sigma X$ is the ``pinch'' map
defined by
\[ [(x,t)]\mapsto \begin{cases} \big([(x,2t)],*\big) & 0<t<1/2 \\
\big(*,[(x,2t-1)]\big)& 1/2\leq t\leq 1 \end{cases} \] and $\nabla\co \Sigma X\vee\Sigma X\to\Sigma X$ is the
``fold'' map. The fact that we can add maps gives the theory a very different, more algebraic flavor, than that of
unstable homotopy theory. Moreover, the same proof as for ordinary homotopy groups shows that $[\Sigma^2 X,Y]$ is
an abelian group.

For us, ``stable'' just means that we can suspend any number of times, even an infinite number, as needed. More
precisely let $QY = \colim \Omega^n \Sigma^n Y$. Then we define the stable homotopy classes
of maps $\{X,Y\} = [X, QY]$. If $X$ is a finite complex $\{X,Y\} = 
\colim{ [\Sigma^n X, \Sigma^n Y]}$. For a discussion of spectra, see Adams \cite{A}.

\subsection{Classifying spaces}
Let $G$ be a finite group. Define $EG$ to be a free, contractible $G$--space, and let $BG$ denote the quotient
$EG/G$. The contractibility of $EG$ shows us that $BG$ is a space with a single nontrivial homotopy group:
$\pi_1(BG)=G$. We give two explicit constructions of $EG$ and then give an application.

\subsubsection{Milnor's definition}
We define \[EG=\bigcup_{n} \underbrace{G*\dots*G}_n,\] where $*$
denotes the join of two spaces, which we take to be the suspension
of the smash product. Since the join includes a suspension, the
greater the number of copies of $G$ being joined, the higher the
connectivity, and so $EG$ is contractible.

As an example, we take $G=\Z/2.$ In this case,
\[E\Z/2 = \bigcup_{n} \underbrace{\Z/2*\dots*\Z/2}_n=\bigcup_{n}
S^{n-1}=S^{\infty},\] and the $\Z/2$ action is the usual diagonal
reflection action, and in this case, the quotient $B\Z/2$ is just
$\mathbb RP^{\infty}.$

\subsubsection{Simplicial model} We can think of a group as a
category with a single object and whose morphisms are the elements
of the group. We can now pull in the categorical construction of
the nerve, and this will give us a model for $BG$.

First we recall briefly the definition of a simplicial set. A more complete reference is May \cite{M}. Let $\Delta$ be
the category whose objects are the sets ${0,\dots,n}$ for all $n$ and whose morphisms are nondecreasing maps. A
\textit{simplicial set} is a contravariant functor from $\Delta$ to the category of sets. We can think of a
simplicial set as a collection of sets indexed by the natural numbers together with a large family of structure
maps called faces and degeneracies which satisfy certain properties, modeled dually on the inclusion of faces in
the standard simplices in $\mathbb R^n$.

To any simplicial set $S$, we can associate a topological space,
the geometric realization, $|S|$. Loosely speaking, this is
defined by putting a copy of the standard $n$--simplex in for every
element of $S_n$ and gluing them all together via the face and
degeneracy maps.

To any category $\mathcal C$, we can associate a simplicial set,
the \textit{nerve}, $N\mathcal C_*$. The $k$--simplices of
$N\mathcal C$ are the $k$--tuples of composable morphisms in
$\mathcal C$. The face maps are induced by the various ways to
compose adjacent maps (or to forget the ends), and the
degeneracies comes from inserting the identity map in various
places. We define the classifying space $B\mathcal C$ to be the
geometric realization of the nerve. With a little work, one can
quickly show that this construction is functorial.

In the case $\mathcal C=G$, a finite group, $N\mathcal C_k=G^k$,
since all morphisms are composable. We then get a model of $BG$ by
taking the geometric realization.

This construction has some very nice advantages over the previous
one, and to show this, we need a small proposition.

\begin{prop}
If $F_0,F_1\co \mathcal C\to\mathcal C',$ and $H$ is a natural
transformation from $F_0$ to $F_1$, then $BF_0\simeq BF_1$ as maps
$B\mathcal C\to B\mathcal C'$, and the homotopy is given by $BH$
on $B(\mathcal C\times \{0\to 1\})=B\mathcal C\times [0,1].$
\end{prop}

This immediately gives us an important result about conjugation.

\begin{cor}\label{c1}
Let $x\in G$, and let $C_x(g)=x^{-1}gx$ denote conjugation by $x$.
We then have $BC_x\simeq \Id_{BG}.$
\end{cor}

\begin{proof}
There is a natural transformation between $C_x,$ viewed as an
endofunctor of $G$, and the identity functor given by
``multiplication by $x$'':
\[ \xymatrix{{e}\ar[r]^{C_x(g)} \ar[d]_{x} & {e}\ar[d]^{x} \\ {e}
\ar[r]_{g} & {e}} \] where $e$ denotes the single object in the
category. In other words, the morphism $x\co e\to e$ is a natural
transformation between $C_x$ and the identity, and the result
follows.
\end{proof}

It is this simple corollary which gives us the basic connection between group theory and the homotopy theory of
classifying spaces.

\subsection{Group cohomology}
The space $EG$ allows us to define group homology and cohomology.
The singular chains $C_*(EG)$ is a $\Z[G]-$free resolution of
$\Z$.

\begin{defn}
For any $G-$module $M$, let $H_*(G;M)=H(C_*(EG)\otimes_{\Z[G]}
M),$ and let $H^*(G;M)=H(\Hom_{\Z[G]}(C_*(EG),M).$
\end{defn}

Note in particular that if $M=\Z$, the trivial $G-$module, then
$H_*(G;M)=H_*(BG)$ and similarly for cohomology. If $M$ is not
trivial, then $H_*(G;M)$ can be similarly related to $H_*(BG)$ but
with twisted coefficients.

In what follows all cohomology is taken with simple coefficients
in $\F_p$.

\section{Stable splittings}
Suppose that $\Sigma BG=X_1\vee\dots \vee X_N$. If we can do this, then for any generalized cohomology theory $E$,
\[ E^*(BG)=E^{*+1}(\Sigma BG)=\bigoplus E^{*+1}X_i.\] In general, this is a simpler object to study. We want now
to find ways to relate the $X_i$ to $G$ itself.

\subsection{Summands via idempotent self-maps}
Let $e\co X\to X$ for some pointed space $X$. If $e^2\simeq e,$ we call $e$ a \textit{homotopy idempotent}.
We now form the mapping telescope $eX=\mathrm{Tel}(X,e)$ which is the homotopy colimit of the diagram $\smash{X\stackrel{e\,}{\vphantom{\cdot}\smash{\rightarrow}}X}$.
More explicitly, we start with the disjoint union $$\coprod_{n \geq 0} (X\times [2n,2n+1])$$ and identify
$(x,2n+1)$ with $(e(x),2n+2)$ and pinch $(*,t)$ to a point.
In this case, \[\pi_*(eX)=\colim\pi_*(X)=e_*\pi_*(X),\] where the
structure maps in the limit are $e_*$. A similar statement holds
for homology.

If $X$ is a suspension, then we can add and subtract maps, and in particular, we can form a map $X\to eX\vee
(1-e)X$ whenever $e$ is a homotopy idempotent. From the above comments, this is an equivalence. Our next task is
then to find idempotents in $[X,X]$. In the case of $X=BG$ or $\Sigma BG$, we shall get the first layer of these
from algebra, using the $\Aut(G)$ action on $[BG,BG]$.

Rather than looking at homotopy classes of maps, we'll look at stable homotopy classes of 
self maps $\{BG,BG\}$. Under composition,  this has the structure of a ring.
The group of stable homotopy self-maps carries an action of $\Aut(G)$ via the map which sends $\alpha\in \Aut(G)$ to the stable class of $B\alpha$. This therefore extends to a map of 
rings from $\Z[\Aut(G)]\to \{BG,BG\}$. If
we can find idempotents in $\Z[\Aut(G)],$ then we can push them forward to stable homotopy 
idempotents. Since we are interested in working one prime at a time and since idempotent
theory is easier for completed rings, we shall assume $BG$ is completed at $p$ and consider 
the induced map $\Z_p[\Aut(G)] \to \{BG, BG\}$. If $G$ is a $p$--group then $BG$ is already 
$p$--complete.

We start by reducing mod $p$, since any idempotent $e\in \F_p[\Aut(G)]$ lifts to an 
idempotent in $\Z\,_p[\Aut(G)].$ Moreover, if we have a primitive orthogonal idempotent 
decomposition $1=e_1+\dots+e_n$ in
$\F_p[\Aut(G)],$ where $e_ie_j=0$ for $i\neq j$, $e_i^2=e_i$ then this lifts to a 
decomposition of the same form in $\Z\,_p[\Aut(G)].$

\begin{exm}
If $G=\Z/p$, then $\Aut(G)=\Z/(p-1)$. If $p=3$, then we can
readily find two idempotents in $\F_3[\Z/2],$ namely $-1-e$ and
$-1+e$, where $e$ is the nontrivial element in $\Z/2.$

In general, there are $p-1$ primitives with idempotents given by
\[ e_i=\prod_{j\neq i} \frac{\xi-a^i}{a^i-a^j}, \quad i=0,\ldots,p-2 \]
where $\xi$ is the generator of $\Z/p-1$ and $a$ is the element in
$\F_p=\Z/p$ by which $\xi$ acts.
\end{exm}

\begin{prop}
Stably and $p$--completed, \[B\Z/p\simeq X_0\vee\dots\vee X_{p-2},\] where $X_i=e_iB\Z/p.$
\end{prop}

We will say more about $X_0$ below. One can also try to use the full
ring of endomorphisms for an abelian $p$--group $P$. This approach has
been thoroughly studied by Harris and Kuhn \cite {HK}.

\section{Transfer maps}

Let $H\subset G$ be a subgroup of index $[G:H]=n.$ If we take the quotient
of $EG$ by $H$, then we get $BH$, since $EG$ is contractible and being
$G$ free forces it to be $H$ free. We can further quotient by all of $G$
to get a map $BH\to BG$ and the fiber of this map is $G/H$. In other
words, we have an $n$--sheeted cover $BH\to BG$. The map $BH\to BG$
is also easily seen to be equivalent to $B$ of the inclusion $H\to G$.

The transfer is a stable map which goes from $BG$ back to $BH$. In
cohomology, we can easily define it. Let $\pi\co X\to Y$ be an
$n$--sheeted cover. For each small enough simplex $\Delta\in C_*(Y)$,
we can find $n$ simplices in $C_*(X)$ lying over it. The transfer is
the map in homology induced by
\vadjust{\vskip-3pt}
\[ \Delta\mapsto\sum_{\Delta'\in\pi^{-1}(\Delta)}\Delta'.
\vadjust{\vskip-3pt}
\]
If we compose now with the projection map, then it is clear that the
composite is simply multiplication by $n.$

Actually getting a stable map requires a little more work. Write
$G=\coprod \tau_i H.$ Given an element $\tau_i$, left multiplication
by $g\in G$ sends it to $\tau_{\sigma(g)(i)} h_{i,g}$. This gives us a
permutation representation $\sigma\co G\to\Sigma_n$ and a homomorphism
\begin{align*}
G&\to H^n\rtimes\Sigma_n=\Sigma_n\wr H \\[-0.5ex]
g&\mapsto \big(h_{1,g},\dots,h_{n,g},\sigma(g)\big).
\end{align*}
We define the \textit{transfer} to be the map adjoint to the
composite
\[ \xymatrix{{BG}\ar[r]\ar[rrrd] & {B(\Sigma_n\wr H)}\ar[r]^{=\quad} &
{BH^n\times_{\Sigma_n} B\Sigma_n}\ar[r] & {(QBH)^n\times_{\Sigma_n} B\Sigma_n} \ar[d]_{\Theta} \\ {} & {} & {} &
{QBH}},\]
where $\Theta$ is the Dyer--Lashof map arising from the infinite
loop structure of $QBH$. It is not difficult to see that in homology
this map agrees with the previous definition for the covering $BH \to
BG$. Actually the map we have defined is sometimes referred to as the
reduced transfer. Let $BG^{+}$ denote $BG$ with an added disjoint
basepoint so that $BG^{+} \simeq BG  \vee S^0$. Then it is easy to
extend this definition to a stable map $\tr \co BG^{+} \to BH^{+}$
which is multiplication by $[G,H]$ on the bottom cell. For a detailed
exposition see Kahn and Priddy \cite{KP}; another approach is given by
Adams \cite{A}.

\subsection{Properties of the transfer and corollaries}
We have already seen homologically that the composite
$BG\xrightarrow{\tr} BH\xrightarrow{Bi} BG$ is multiplication by
the index $[G:H].$

\begin{cor}\label{c2}
If $H\subset G$, and $[G:H]$ is prime to $p$, then $BG$ is a
stable summand of $BH$ when completed at $p$.
\end{cor}
\begin{proof}
Since $[G:H]$ is prime to $p$, it is a unit in $\Z\,_p$, and
multiplication by it is an equivalence. Thus the transfer and
inclusion give the splitting.
\end{proof}
\begin{cor}
Stably and completed at $p$, $B\Sigma_p$ is a stable summand of $B\Z/p$.
\end{cor}

\begin{prop}[Properties of the transfer] We will write $\tr_H^G$
for the transfer with $H$ considered as a subgroup of $G$.
\begin{enumerate}
\item If $H=G$, then $\tr_H^G=\Id$.

\item If $K\subset H\subset G$, then $\tr_K^G=\tr_H^G\circ \tr_K^H.$

\item The transfer is natural with respect to maps of coverings.

\item The ``Double Coset Formula'' holds: If $H,K\subset G$, write $G=\coprod KxH$ for some collection of $x\in
G$. Let $\tr_x$ denote the transfer $BK\to B(K\cap x^{-1}Hx),$ and let $i_x$ denote the inclusion $xKx^{-1}\cap
H\to H$. Then if $i_K$ is the map $BK\to BG$, we have
\[ \tr_H^G\circ i_K=\sum_x i_x\circ C_{x^{-1}}\circ \tr_x. \]
\end{enumerate}
\end{prop}

\begin{lem}
If $G$ is an elementary abelian $p$--group, and $H\subsetneq G$, then the transfer induces the zero map in mod $p$
cohomology.
\end{lem}
\begin{proof}
The map $i_H^*$ is surjective, since $H$ sits inside $G$ as a
summand. Since \[{\tr_H^G}^*i_H^*(x)=[G:H]x=p^nx=0,\] we conclude
that $\tr^*=0.$
\end{proof}
\begin{cor}\label{c4}
If $V\subset G$ is elementary abelian, then \[i_V^*\circ
{\tr_V^G}^*=\sum_{w\in N(V)/V} C_{w}^*.\]
\end{cor}
\begin{proof}
Let $K{=}V$ in the double coset formula and vary $x$ over coset
representatives.\!\!
\end{proof}

It follows from \fullref{c2} that if $G_p$ is a Sylow $p$--subgroup, then $BG$ is a stable summand of $BG_p$ after
$p$--completion. We now specialize to the case that $G_p=V$ is elementary abelian.

\begin{thm}\label{t1}
If $V\subset G$ is an elementary abelian Sylow $p$--subgroup, then
\begin{enumerate}
\item $H^*(G)\cong H^*(V)^W,$ where $W=N_G(V)/V$ is the Weyl group. \item $BN_G(V)\to BG$ is an
$H\Z/p$--equivalence, even unstably.
\end{enumerate}
\end{thm}
\begin{proof}
The second result follows immediately from the first.

For the first part, note that we always have $H^*(G)\subset
H^*(V).$ Since conjugation acts as the identity on $H^*(G)$, we
must have \[H^*(G)\subset H^*(V)^{N_G(V)}.\] Since conjugation by
$V$ is trivial in cohomology, $H^*(V)^{N_G(V)}=H^*(V)^{W}.$

From \fullref{c4}, we know that the composite
\[H^*(V)\xrightarrow{\tr^*} H^*(G)\xrightarrow{i^*} H^*(V)^W\] is
just $\sum_{w\in W} C_w.$ Since $|W|$ is prime to $p$, it is invertible in $\Z_p,$ and $e=\sum C_w/|W|$ is an
idempotent invariant under the action of $W$. Conversely, all invariants arise in this way, since on the
subalgebra of $W$--invariants, the composite is just multiplication by $|W|$ and is therefore invertible.  This in
particular shows that $i^*$ is surjective, and the result follows.
\end{proof}

\begin{exm}
For $G= \Sigma_p$, the Sylow $p$--subgroups are $\Z/p$, and $N_{\Sigma_p}(C_p)/C_p=\Z/(p-1).$ Now
$H^*(\Sigma_p)=H^*(\Z/p)^W$. The group $W$ acts on the cohomology $H^*(\Z/p)=E(x_1)\otimes\F_p[y_2]$ as multiplication by a generator of
$\F_p^\times$ on $x_1$ and $y_2=\beta(x_1).$ The fixed point algebra is then generated by $\smash{x_1y_2^{p-2}}$ and
$\smash{y_2^{p-1}}$ as an unstable algebra over the Steenrod algebra. This shows that \[H^*(\Sigma_p;\F_p)=\begin{cases}
\Z/p & *=0,-1\!\mod\! p
\\ 0 & \text{otherwise.} \end{cases} \] These dimensions explain
why the Steenrod operations occur where they do, just as a similar computation for the map $B(\Z/p \times \Z/p)
\to B(\Z/p \wr \Z/p) \to B\Sigma_{p^2}$ yields the Adem relations.
\end{exm}

With more work one can show a generalization of \fullref{t1}.

\begin{thm}[Harris--Kuhn \cite{HK}]
If $P\subset G$ is a Sylow $p$--subgroup, and $P$ is an abelian
$p$--group, then
\begin{enumerate}
\item $H^*(G)\cong H^*(P)^W=H^*(N_G(P))$.

\item $BN_G(P)\to BG$ is an $H\F_p$--equivalence.
\end{enumerate}
\end{thm}

\subsection{Modular representation theory}
If $p$ divides the order of the automorphism group, then the representation theory of $\Aut(G)$ over $\F_p$ lies
in the realm of modular representation theory and hence becomes more complicated. We demonstrate this with some
basic examples of increasing trickiness.

\begin{exm}\label{e3}
Let $G = V_2=\Z/2\times\Z/2$, and take $p=2$. Now
\[\Aut(V_2)=GL_2(\F_2)=\Sigma_3\] has order divisible by $2$, and we can
find simple generators \[ \sigma=\left[\begin{matrix} 0 & 1
\\ 1 & 1 \end{matrix}\right] \quad\text{and}\quad \tau=\left[\begin{matrix}
0 & 1 \\ 1 & 0 \end{matrix}\right]. \] In $R=\F_2[\Aut(V_2)]$ the
element
\[f_1=1+\sigma+\sigma^3\] is an idempotent, since $f_1^2=3f_1$.
Since $\tau\sigma\tau=\sigma^2$, $f_1$ and $f_2=1-f_1$ are central
idempotents. With a small bit of work, one can show the following.
\begin{prop}
$R\cong Rf_1\times Rf_2$, and $Rf_1=E(\gamma)$, where
$\gamma=\sum_{g\in GL} g$ and $Rf_2=M_2(\F_2)$.
\end{prop}

We can lift $f_1$ and $f_2$ to idempotents $e_1$ and $e_2$ in $\{BV_2,BV_2\}$, so we conclude that \[
BV_2=e_1BV_2\vee e_2BV_2.\] The first summand we can identify, as it is clearly the same as $B$ of the semi-direct
product \[(\Z/2\times\Z/2)\rtimes \Z/3=A_4.\]
Now $f_2$ can be written as the sum of two idempotents $F_1$ and
$F_2$, where
\[F_1=(1+\tau\sigma)(1+\tau)\quad\text{and}\quad
F_2=(1+\tau)(1+\tau\sigma).\] These can also be lifted to idempotents in $\Z\,_2[\Aut(V)]$, so $BV_2$ splits
further as
\begin{equation}
\label{eq:1}  BV_2=BA_4\vee\tilde{e}_1BV_2\vee\tilde{e}_2BV_2.
\end{equation}
Finally, $\tilde{e}_1BV_2 \simeq \tilde{e}_2BV_2,$ since we have a
sequence
\[F_1R\xrightarrow{1+\tau} F_2R\xrightarrow{1+\tau\sigma}
F_1R\xrightarrow{1+\tau} F_2R\] in which the composite of any two successive arrows is the identity. It is known
that
\[\tilde{e}_1BV_2=L(2)\vee B\Z/2,\] where $L(2)$ is a spectrum that
is related to Steenrod operations of length two. We therefore have the following result of Mitchell \cite{MiP}:
\[BV_2=BA_4\vee (L(2)\vee\mathbb RP^{\infty})\vee (L(2)\vee\mathbb
RP^{\infty}),\]in which the summands are indecomposable.

\end{exm}

\begin{exm}\label{e4}
Let $G=D_8$, the dihedral group of order $8$. One can show that $\Aut(D_8)=D_8$, so this is a two group.
\begin{lem}
If $G$ is a $p$--group, then $\F_p[G]$ has only one simple module, the trivial one.
\end{lem}
The lemma follows from the fact that the augmentation ideal is nilpotent in this case.

In the case $G=D_8$, the lemma shows that we have only one idempotent in $\F_2[\Aut(G)]$, namely the element $1$.
Nevertheless, we have the following splitting \cite{MiP}:
\[BD_8=BPSL_2(\F_7)\vee(L(2)\vee\mathbb RP^{\infty})\vee (L(2)\vee\mathbb
RP^{\infty}).\]
\end{exm}

Using ring theory, we can get a more direct relationship between the structure of $R=\F_p[\Aut(G)]$ and
idempotents. Let $J$ be the Jacobson radical of $R$, namely the elements annihilating all simple $R-$modules or
equivalently the intersection of all maximal ideals. Ring theory tells us that $R/J$ is semisimple and therefore
splits as a product of matrix rings over division algebras:
\[R/J=M_{n_1}(D_1)\times\dots\times M_{n_k}(D_k).\] The simple
$R$--modules are just the columns of the various matrix rings. Lifting the idempotents from this decomposition, we
obtain in $R$ a primitive orthogonal decomposition \[1=\sum_j e_j=n_1e_1+\dots+n_ke_k.\]

For the earlier \fullref{e3} of $V_2$, the Jacobson radical is the
ideal generated by the element $\gamma$, and so
\[\F_2[\Aut(V_2)]/J=\F_2\times M_2(\F_2),\] giving us the correct
number of factors for the decomposition of $BV_2$ in \eqref{eq:1}. To obtain the complete decomposition of
 \eqref{eq:1} one must use the full endomorphism group $\End(V_2)$.

\begin{exm}
Next, we look at a more complicated example, $V_3=(\Z/2)^3.$ In this case, $R=\F_2[GL_3(\F_2)]$ has $4$ simple
modules.

\begin{center}
\begin{tabular}{|c|c|c|c|c|}
\hline Module & $\F_2$ & $V_3$ & $V_{\smash[t]{3}}^*$ & $\St_3$
\\ \hline Dimension & $1$ & $3$ & $3$ & $8$ \\ \hline
\end{tabular}
\end{center}
where ${V_3}^*$ is the contragredient module and $\St_3$ is the Steinberg module (described explicitly
in \fullref{e6} below). This means that we have a splitting
\[BV_3=(\tilde{e}_1BV_3)\vee 3(\tilde{e}_2BV_3)\vee
3(\tilde{e}_3BV_3)\vee 8(\tilde{e}_4BV_3).\] As before, we can
identify the first summand. In this case, it is the same as $B$ of
the group $\F_8^\times\rtimes\, \mathrm{Gal}(\F_8/\F_2)$ obtained by taking the 
semidirect product of $V_3=\F_8$ with $\Z/7\rtimes\Z/3=\mathrm{Gal}(\F_8/\F_2)$.
\end{exm}

Generalizing the example of the dihedral group $D_8$ (\fullref{e4}), we consider
$P=D_{2^n}$, the dihedral group of order $2^n$. We have the following
splitting theorem:
\begin{thm}[Mitchell--Priddy \cite{MiP}]
$BD_{2^n}=BPSL_2(\F_q)\vee 2L(2)\vee 2B\Z/2$, where $q=p^k$ for $p$ odd 
and $D_{2^n}$ is a Sylow $2$--subgroup of
$PSL_2(\F_q)$.
\end{thm}
This condition translates to saying that $n=\nu_2(\tfrac{q^2-1}{2})$ is the order of $2$ in $(q^2-1)/2$.

We can also explain the existence of the summands $2L(2)\vee
2B\Z/2$. There are two nonconjugate copies of $\Z/2\times\Z/2$
sitting inside $D_{2^n}$, and the summands in question appear via
transfers from $BD_{2^n}$ to $B$ of these subgroups, followed by
projection onto the summands.

\section{The Segal Conjecture and its consequences}

While studying the idempotents in the ring $\{BG,BG\}$ via $\Aut(G)$ provides a good bit of 
information about the splittings of $BG$, if $G=P$ is a $p$--group then the Segal conjecture 
completely determines the ring $\{BG,BG\},$ so it is to this that we turn.

Let $A(P,P)$ be the Grothendieck ring of $P\times P$ sets which are free
on the right. The sum is given by disjoint union, and the product is the
product over $P$.  If $P_0\subset P$ is a subgroup, and $\rho\co P_0\to P$
is a homomorphism, then we can define elements of $A(P,P)$ by
\[ P\times_\rho P=P\times P/ \sim \]
where $(xp_0,y)\sim (x,\rho(p_0)y)$, $p_0 \in P_0$. As a group, $A(P,P)$ is a 
free abelian group on these transitive sets. There is a homomorphism 
\begin{gather*}
\alpha\co A(P,P) \to \{BP, BP\} \\
\alpha(P\times_\rho P) =  (BP_+\xrightarrow{\tr_{P_0}^P}
{BP_0}_+\xrightarrow{B\rho} BP_+).\qquad\qquad\tag*{\hbox{defined by}}
\end{gather*}
Upon completion this map is essentially an isomorphism.
More precisely $A(P,P)$ contains an ideal $\tilde{A}(P,P)$ which is a free abelian group on 
the classes $P\times_\rho P - (P/P_0 \times P)$ and we have the following theorem.

\begin{thm}[Carlsson \cite{C}; Lewis--May--McClure \cite{LMM}]
The map $\alpha$ induces a ring isomorphism
$$
\tilde{\alpha}\co  \tilde{A}(P,P)\otimes \Z_p \to \{BP,BP\}. 
$$
\end{thm}

\begin{cor}\label{c5}
\begin{enumerate}\item $\{BP,BP\}$ is a finitely generated, free $\Z_p$--module. 
\item $BP$ splits as a finite wedge $X_1\vee\dots
X_n$, where the $X_i$ are indecomposable $p$--complete spectra, unique up to order 
and equivalence.
\end{enumerate}
\end{cor}
\begin{proof}
The first is immediate. For the second, after tensoring with $\F_p$ we have a finite dimensional $\F_p$--algebra.
This means that $1=\sum e_i$ is a decomposition into primitive idempotents unique up to order and conjugation.
\end{proof}

\begin{cor}
Given a finite $p$--group $P$, there exist finitely many stable homotopy types of $BG$ with $P$ a Sylow
$p$--subgroup of $G$.
\end{cor}
\begin{proof}
We have already seen that if $P$ is a Sylow $p$--subgroup of $G$, then we have a splitting $BP\simeq
BG\vee\text{Rest}$. The finiteness result of the \fullref{c5} gives this one.
\end{proof}

\begin{cor}
Each summand of $BP$ is also an infinite complex.
\end{cor}
\begin{proof}
By the first part of \fullref{c1}, we know that $\{BP,BP\}$ is torsion free. If $X$ were both a finite
complex and a summand of $BP$, then $BP\to X\to BP$, the projection followed by the inclusion, would be a torsion
free map. However, if $X$ is a finite complex, then the identity map of $X$ has torsion, since $X$ is
$p$--complete, so the torsion free composite must as well.
\end{proof}

\begin{thm}[Nishida \cite{N}]\label{t5}
Given $G, G'$ finite groups with $BG\simeq BG'$ stably at $p$, then the Sylow $p$--subgroups of $G$ and $G'$ are
isomorphic.
\end{thm}
We shall derive this from a more general result in \fullref{s6}.

\subsection{Analysis of indecomposable summands}
We will now assume that $X$ is an indecomposable summand of $BP$ for $P$ a fixed $p$--group.

\begin{defn} $X$ \textit{originates} in $BP$ if it does not occur
as a summand of $BQ$ for any $Q\subsetneq P$. $X$ is a \textit{dominant summand} of $BP$ if it originates in $BP$.
\end{defn}

The notion of dominant summand is due to Nishida. As an example, for $BV_2$ the dominant summands are $BA_4$ and
$L(2)$.

It is also clear that every $X$ must originate in some subgroup $Q$ of $P$.

Now let $J(P)$ be the ideal in $\{BP,BP\}$ generated by the maps which factor through 
$BQ$ for some $p$-group $Q$ such that $|Q| < |P|$.
In other words, these are the maps which arise from transitive sets for which $P_0$ is a proper subgroup or if $P_0=P$ from proper 
(ie nonsurjective) endomorphisms
$\rho\co P_0 \to P$. Since every summand $X=eBP$ for some idempotent $e\in\{BP,BP\}$, 
$X$ is dominant if and only if $e\notin J(P)$. Furthermore we have an isomorphism of rings
\[\Z_p[\Out(P)]\xrightarrow{i}\{BP,BP\} \xrightarrow{\pi} \{BP,BP\}/J(P)=\Z_p[\Out(P)].\]
This follows by remembering that $\{BP,BP\}$ is generated by transfers followed by homomorphisms. If the subgroup
for the transfer is proper, or if the homomorphism is not surjective then this 
map is in $J(P)$, so all that we have left over are the automorphisms of
$P$. This gives us the following equivalence:
\begin{multline*}\{\text{Homotopy types of dominant summands of $BP$}\}\rightarrow\\
\{\text{Isomorphism classes of simple
$\F_p[\Out(P)]$--modules}\}\end{multline*}
given by $eBP\mapsto e_0\F_p[\Out(Q)]$, where $e_0$ is determined as follows: since $e$ is primitive we can find a
primitive idempotent $\tilde{e} \in \Z_p[\Out(P)]$ such that $i(\tilde{e}) = e$ mod $J$. Then $e_0$ is the mod $p$
reduction of $\tilde{e}$.

Let $S$ be a simple $\F_p[\Out(Q)]$ module and let $X_{Q,S}$ be its corresponding dominant summand of $BQ$.
\begin{thm}[Mitchell--Priddy \cite{MP1}; Benson--Feshbach \cite{BF}]\label{t6}
There is a complete stable decomposition unique up to order and equivalence of factors
$$BP = {\bigvee\!}_{Q,S}\ n_{Q,S} X_{Q,S}$$
where $Q$ runs over the subgroups of $P$, $S$ runs over the simple $\F_p[\Out(Q)]$ modules, and $n_{Q,S}$ is the
multiplicity of $X_{Q,S}$ in $BP$.
\end{thm}

\subsection{Principal dominant summand}
Among all dominant summands, there is a distinguished one corresponding to the trivial module. If we decompose
$1\in\F_p[\Aut(P)]$ into primitive orthogonal idempotents, then we can consider the image of them under the
augmentation ring map $\epsilon \co \F_p[\Aut(P)]\to F_p$ defined by sending all $g\in \Aut(P)$ to $1\in\F_p$. Since
$e_i^2=e_i$, these must map to either $0$ or $1$ under the augmentation map. Additionally, exactly one must map to
$1$, since the augmentation sends $1$ to $1$, and $e_ie_j=0$ for $i\neq j$. We denote by $e_0$ the idempotent that
maps to $1$ and say that it is the \textit{principal} idempotent. The corresponding summand will be denoted $X_0$
and called the \textit{principal dominant summand}.

\begin{prop}[Nishida \cite{N}]
$X_0$ is a summand of $BG$ for all $G$ with $P$ a Sylow $p$--subgroup.
\end{prop}

For two examples, for $B\Z/p$, the principal dominant summand is $B\Sigma_p$, and for $BV_2$, it is $BA_4$.

\subsection{Ring of universally stable elements}

For a fixed $p$--group $P$, we define the ring of universally stable elements as
\[I(P)=\bigcap_{\substack{G\supset P \text{as}\\\text{a Sylow}\\\text{$p$--subgroup}}} \mathrm{Im}\big(H^*(G)\to
H^*(P)\big).\]
\begin{thm}[Evens--Priddy \cite{EP}]\label{t7}
$H^*(P)$ is a finite module over $I(P)$. This implies that $H^*(P)$ is a finitely generated algebra over $I(P)$ of
the same Krull dimension, which in turn implies $I(P)$ is a finitely generated $\F_p$--algebra.
\end{thm}

This follows from Quillen's theorem \cite{Q}.

\begin{prop}{\rm \cite{EP}}\qua
If $E$ is an elementary abelian $p$--group of rank $n$, then
$I(E)=H^*(E)^{GL(E)},$ except when $p=n=2$, where it is
$H^*(BA_4).$
\end{prop}

In general, $I(P)$ is not realizable as the cohomology of a spectrum. In certain 
familiar cases, however, it is
not only realizable but also connected to the cohomology of the principal dominant summand.

\begin{prop}{\rm \cite{EP}}\qua 
We have the following rings of universally stable
elements.
\[\begin{array}{lll}
I(D_{2^n})=H^*(PSL_2(\F_q)) &&
n=\nu_2\left(\tfrac{q^2-1}{2}\right) \\
I(Q_{2^{n+1}})=H^*(SL_2(\F_q)) && n=\nu_2\left(\tfrac{q^2-1}{2}\right) \\
I(SD_{2^n})=H^*(SL_3(\F_q)) && n=\nu_2((q^2-1)(q+1))\vphantom{\tfrac{q^2-1}{2}},\ \ q=3\!\mod
4
\end{array}\]
 For all of these groups, $X_0$ is $B$ of the group shown.
\end{prop}

\section{Summands of supergroups}
Recall that if $X$ is an indecomposable summand of $BQ$, then we
say it is dominant if it is not a summand of $BQ'$ for any
subgroup $Q'\subsetneq Q$. If $P$ is a $p$--group, and $Q\subset P$ is
a subgroup with $X$ a dominant summand of $BQ$, we can ask when
$X$ is a summand of $BP$.

\begin{exm}\label{e6}
We recall the definition of the Steinberg module. Define ${e_{\St_n}}\in R = \F_p[GL_n(\F_p)]$ by
\[ e_{\St_n} = \frac{1}{[GL_n(\F_p):U_n]}\bigg(\sum_{b\in B_n, \sigma\in\Sigma_n} (-1)^\sigma
b\sigma\bigg),\] where $B_n$ is the Borel subgroup of upper triangular matrices and $U_n$ is the unipotent
subgroup thereof. This element is idempotent and primitive, and if we let $\smash{\St_n=e_{\St_n}R}$, then $\St_n$ is a
simple, projective module of dimension $\smash{p^{\binom{n}{2}}}$.

If $V$ is an elementary abelian group of dimension $n$, then
\[e_{\St_n}BV=L(n)\vee L(n-1),\] where
$L(n)=\Sigma^{-n}\big(Sp^{p^n}(S^0)/Sp^{p^{n-1}}(S^0)\big)$ and $Sp^n(X)= X^n/\Sigma_n$
is the symmetric product
\cite{MiP}.
\end{exm}

\begin{prop}{\rm\cite{P}}\qua
If $P$ is a $p$--subgroup of rank $n$ (ie $n$ is largest number such that $(\Z/p)^n$ is a subgroup of $P$), then
$L(n)$ is a summand of $BP$ if and only if $P$ contains a self-centralizing elementary abelian subgroup of rank $n$.
\end{prop}

\begin{exm} For $p=2$ and low values of $n$, we have the following
examples.

$n=2$\qua $L(2)$ is a summand of $BD_{2^n}$. 

$n=1$\qua $L(1)$ is actually $\mathbb RP^{\infty}$, and we
have seen already instances when this occurs as a summand. 

$n=1$\qua As a ``non-example'', $L(1)$ is not a
summand of $BQ_8$, since the $\Z/2\subset Q_8$ is central.
\end{exm}

Now recall the principal summand $X_0$. Let $PS(X_0,t)$ be its Poincar\'{e} series in mod--$p$ cohomology. By the
dimension of $X_0$ we mean the order of the pole of $PS(X_0,t)$ at $t=0$.

\begin{prop}[Martino--Priddy \cite{MP3}]
$X_0$ has dimension $n$ in $BP$ if the rank of $P$ is $n$.
\end{prop}

\subsection{Generalized Dickson invariants}
For the Dickson invariants, we normally start with the algebra $H^*(E)^{GL(E)}$. If $E$ is an elementary abelian
$p$--group of rank $n$, then we can form the composite
$$E\xrightarrow{\Reg}\Sigma_E\to U(p^n),$$ where $\Reg$ is the
regular representation of $E$ acting on itself. In cohomology, the Chern classes $\smash{c_{p^n-p^i}}$ map to
the  Dickson invariants $c_i$. These Dickson invariants carve out a polynomial invariant subalgebra of
$\smash{H^*(BE)^{GL(E)}}$ of dimension $n$.

For a general $p$--group $P$, we can formally mirror the above construction. Let $\rho$ denote the composite
complex representation
\[\rho\co  P\xrightarrow{\Reg}\Sigma_{|P|}\to U(|P|),\] and define the
 (generalized) Dickson invariants of $P$ to be the image under $\rho^*$ of
\[\F_p[c_{p^s(p^n-p^i)}\,:\,i=0,\dots, n-1\,;\, s=|P|/p^n]\subset
H^*(U(|P|)\] where $s=|P|/p^n$.
This obviously forms a subalgebra of $H^*(BP)$ which
we will denote $D(P)$.

\begin{prop} If $P$ is a $p$--group of rank $n$, then
\begin{enumerate}
\item $D(P)$ is a polynomial ring of dimension $n$. \item $D(P)\subset H^*(BP)^{\Out(P)}$.
\end{enumerate}
\end{prop}

\begin{proof}
This is easy to see. For the first part, we use the fact that the
composite of the inclusion of $E$, an elementary abelian subgroup
of rank $n,$ into $P$ followed by the regular representation map
to $\Sigma_{|P|}$ is the same as $p^s$ times the regular
representation map of $E$. By the naturality of cohomology, the
classes $c_{p^s(p^n-p^i)}$ pull back to the classes
$\smash{(c_{p^n-p^i})^{p^s}}$.

For the second part, given an automorphism $f$ of $P$, we can form
a commutative square
\[\xymatrix{{P}\ar[d]_f \ar[r] & {\Sigma_{P}}\ar[d]^{C_f} \\ {P}
\ar[r] & {\Sigma_P}}\] As before conjugation induces the identity in the cohomology of the classifying space;
this shows that the pullback of any classes coming from $H^*(\Sigma_P)$ lies in the invariants of $H^*(P)$ under
$\Out(P).$
\end{proof}

\section[Stable classifications of BG at p]{Stable classifications of $BG$ at $p$}
\label{s6}
The following result gives a classification of the stable type of $BG$ in terms of its 
$p$--subgroups $Q$ and
associated $\Out(Q)$ modules. Let $\Rep(Q,G)=\Hom(Q,G)/G$ and $\Inj(Q,G)\subset \Rep(Q,G)$be the 
classes of injections. Let $\Cen(Q,G)\subset \Inj(Q,G)$ be represented by monomorphisms 
$\alpha\co Q \to G$ such that $C_G(\mathrm{Im}\alpha)/Z(\mathrm{Im}\alpha)$ is a $p'$--group.

\begin{thm}[Martino--Priddy \cite{MP2}]
Let $G$ and $G'$ be finite groups. The following are equivalent:
\begin{enumerate}
\item $BG\simeq BG'$ stably at $p$. 

\item For every finite $p$--group $Q$, there is an isomorphism of $\Out(Q)$--modules
\[\F_p[\Cen(Q,G)]\cong\F_p[\Cen(Q,G')].\]
\end{enumerate}
\end{thm}

{\bf Note}\qua The proof of a related classification result of \cite{MP2} contains 
an error. See \cite{MP4} for a correction.
\medskip

A stable equivalence $BG \simeq BG'$ at $p$ induces an isomorphism 
$$  \F_p[\Inj(Q,G)] \cong \F_p[\Inj(Q,G')] $$
of $\Out(Q)$ modules. From this we can easily derive Nishida's result, \fullref{t5}:
\begin{cor} If $BG \simeq BG'$ stably at $p$, then $P \cong P'$ where $P$, $P'$ are respective
Sylow $p$--subgroups.
\begin{proof}
 We have
\[0\neq\F_p[\Inj(P,G)]\cong\F_p[\Inj(P,G')].\] This implies $P\subset P'$. Reversing $P$ and $P'$
gives the desired conclusion $P \cong P'.$
\end{proof}
\end{cor}

\begin{cor}
If $BG\simeq BG'$ stably at $p$, then $G$ and $G'$ have the same number of conjugacy classes of $p$--subgroups of
order $|Q|$ for all $Q$.
\end{cor}

\begin{proof}
It is easy to see that
\[\F_p[\Inj(Q,G)]=\bigoplus_{\substack{(Q_1),\ Q_1\cong
Q\\ Q_1\subset G}}\F_p[\Out(Q)]\otimes_{\F_p[W_G(Q_1)]}\F_p,\]
\[W_G(Q_1)=N_G(Q_1)/C_G(Q_1)\cdot Q_1\leqno{\hbox{where}}\] is the ``Weyl group'' of
$Q_1$. From the third part of the theorem, we have an isomorphism
\[\bigoplus_{\substack{(Q_1),\ Q_1\cong
Q\\ Q_1\subset G}}\F_p[\Out(Q)]\otimes_{\F_p[W_G(Q_1)]}\F_p\cong \bigoplus_{\substack{(Q_1),\ Q_1\cong Q\\ Q_1\subset
G'}}\F_p[\Out(Q)]\otimes_{\F_p[W_G(Q_1)]}\F_p,\] and if we apply to this the functor
$\F_p\otimes_{\F_p[\Out(Q)]}(\cdot)$, then we conclude that \[\bigoplus_{\substack{(Q_1),\ Q_1\subset G\\ Q_1\cong
Q}}\F_p\cong\bigoplus_{\substack{(Q_1),\ Q_1\subset G'\\ Q_1\cong Q}}\F_p.\proved\]
\end{proof}

\begin{defn}
Let $H,K$ be subgroups of $G$. We say that $H$ and $K$ are
\emph{pointwise conjugate} in $G$ if there is a bijection of sets
$\smash{H\stackrel{\phi}{\vphantom{-}\smash{\rightarrow}}K}$ such that $\phi(h)=g(h)hg(h)^{-1}$ for some $g(h)\in G$ 
depending on $h$. This is
equivalent to the statement that \[|H\cap (g)|=|K\cap (g)|\] for
all $g\in G$.
\end{defn}

\begin{cor}
Assume that $G$ and $G'$ have normal Sylow $p$--subgroups $P$ and $P'$ respectively. Then $BG\simeq BG'$ stably at
$p$ if and only if there is an isomorphism $\smash{P\stackrel{\phi}{\vphantom{-}\smash{\rightarrow}}P'}$ such that $W_G(P)$ is pointwise conjugate to
$\phi^{-1}W_G(P')\phi.$
\end{cor}

\bibliographystyle{gtart}
\bibliography{link}

\end{document}